\documentclass[10pt]{amsart}
\usepackage{amsfonts}
\usepackage{color}
\usepackage[colorlinks,
            linkcolor=red,
            anchorcolor=blue,
            citecolor=green
            ]{hyperref}

\definecolor{c20}{rgb}{0.,0.7,0.}
\definecolor{c30}{rgb}{0.,0.,1.}
\definecolor{c40}{rgb}{1,0.1,0.7}
\definecolor{c50}{rgb}{1,0,0}
\definecolor{c60}{rgb}{1,0.9,0.1}

\def\ehe#1{\textcolor{c30}{#1}}
\def\ehe#1{#1}

\newcommand{\abs}[1]{\left\lvert #1 \right\rvert}

\newcommand{\E}[1]{\mathbb{E}\left\{#1\right\}}

\newcommand{\pk}[1]{\mathbb{P} \left \{#1 \right \} }

\newcommand{\R}{\mathbb{R}}

\newcommand{\BQN}{\begin{eqnarray}}
\newcommand{\EQN}{\end{eqnarray}}
\newcommand{\BQNY}{\begin{eqnarray*}}
\newcommand{\EQNY}{\end{eqnarray*}}

\newcommand{\BS}{\begin{sat}}
\newcommand{\ES}{\end{sat}}
\newcommand{\BT}{\begin{theo}}
\newcommand{\ET}{\end{theo}}
\newcommand{\BL}{\begin{lem}}
\newcommand{\EL}{\end{lem}}
\newcommand{\BK}{\begin{korr}}
\newcommand{\EK}{\end{korr}}

\newcommand{\BD}{\begin{de}}
\newcommand{\ED}{\end{de}}

\newcommand{\BIT}{\begin{itemize}}
\newcommand{\EIT}{\end{itemize}}
\newcommand{\BDI}{\begin{description}}
\newcommand{\EDI}{\end{description}}

\newcommand{\BRM}{\begin{remarks}}
\newcommand{\ERM}{\end{remarks}}

\newcommand{\BEX}{\begin{example}}
\newcommand{\EEX}{\end{example}}

\newtheorem{theo}{Theorem}[section]
\newtheorem{sat}[theo]{Proposition}
\newtheorem{de}[theo]{Definition}
\newtheorem{lem}[theo]{Lemma}

\newtheorem{example}[theo]{Example}
\newtheorem{korr}[theo]{Corollary}

\newtheorem{remarks}[theo]{Remarks}

\newcommand{\prooftheo}[1]{ \textsc{\bf Proof of Theorem} \ref{#1}:}

\newcommand{\prooflem}[1]{\textsc{\bf Proof of Lemma} \ref{#1}:}
\newcommand{\proofkorr}[1]{\textsc{\bf Proof of Corollary} \ref{#1}:}

\newcommand{\COM}[1]{}

\newcommand{\QED}{\hfill $\Box$}

\date{}

\def\P{ \mathbb{P}}

\def\E{\mathbb{E}}

\def\equivdis{\stackrel{d}{=}}

\begin{document}

\title[The Boundary Non-Crossing Probabilities For Slepian Process]{\ehe{The Boundary Non-Crossing Probabilities for Slepian process}}

\author{Pingjin Deng}
\address{Pingjin Deng, School of Finance, Nankai University, 300350, Tianjin, PR China,
and Department of Actuarial Science,
University of Lausanne\\
UNIL-Dorigny, 1015 Lausanne, Switzerland
}
\email{Pingjin.Deng@unil.ch}

 \maketitle
 {\bf Abstract}:  In this contribution we derive an explicit formula for the boundary non-crossing probabilities for Slepian processes associated with the piecewise linear boundary function. This formula is used to develop an approximation formula to the boundary non-crossing probabilities for general continuous boundaries. The formulas we developed are easy to implement in calculation the boundary non-crossing probabilities.

\bigskip

\date{\today}
\section{ Introduction} \label{Introduction}
A Slepian process $S_a(t),t\in [0,T]$ is a centered, stationary Gaussian processes with continuous paths and covariance function given by
\BQNY
 \E{(S_a(t_1)S_a(t_2))}=
\begin{cases}
 1-\abs{\frac{t_1-t_2}{a}}, & \mbox{if } \abs{t_1-t_2}<1 \\
  0, & \mbox{otherwise},
\end{cases}
\EQNY
with $a$ a positive constant. As mentioned in \cite{slepian1961} and \cite{shepp1971first}, we can represent $S_a$
in terms of standard Brownian motion as the moving window process, namely
\BQNY
S_a(t)=\frac{1}{\sqrt{a}} (B(t+a)-B(t)),\quad t\in [0,T].
\EQNY
Since for any $a>0$ fixed, we have the following
\BQNY
\frac{1}{\sqrt{a}} (B(t+a)-B(t))\equivdis B(\frac{t}{a}+1)-B(\frac{t}{a}),\quad t\in [0,T],
\EQNY
where $\equivdis$ means equivalence in distribution, then denoting  $t'=\frac{t}{a}, T'=\frac{T}{a}$ we obtain the following representation
 $$\{S_a(t);t\in[0,T]\}\equivdis \{S_1(t');t'\in[0,T']\}.$$
Therefore, the distributional properties of the Slepian process $S_a,a>0$ can be derived by those of $S_1$. In the literatures, we always set $T'=1$, and denote
$$S(t):=S_1(t)=B(t+1)-B(t), \quad t\in [0,1].$$
\newline
The process $S(t),t\in [0,1]$ has been extensively studied in several different areas ever since it was defined by Slepian in \cite{slepian1961}. For instance, the small ball probability of $S$  is discussed in  \cite{nikitin2006exact} and \cite{gao2007small}, while the Karhunen-Loève expansion for this process was derived independently in \cite{jin2014gaussian} and \cite{liu2014karhunen}. From the application point of view, Cressie \cite{cressie1980asymptotic} pointed out the importance of $S$ in scan statistics and gave a concrete application of this process in the above mentioned reference. Recently,  Bischoff and Gegg demonstrated that $S$ is useful in the signals detection problems, see \cite{Bischoff2016}.

The Boundary non-crossing probabilities problem, which is studied extensively both in stochastic processes and statistics, see e.g., \cite{MR2016767, MR2028621, MR2175400, MR2591910, MR2028621, hashorva2015boundary} is of interest (and tractable) in the particular case of the Slepian process $S$.  Specifically, in this contribution we are concerned with the boundary non-crossing probability
\BQN\label{EQ1}
p_f:= \pk{S(t)\leq f(t),\;\text{for all}\; t\in [0,1]},
\EQN
where $f$ is a given deterministic measurble function.\\
The computation of $p_f$ is of interest in various statistical applications. For example, when we detect the signal in radars, the test statics we use turn out to be the non-crossing probabilities of Slepian process with nonlinear boundaries, see e.g., \cite{zakai1969threshold}. Though many methods and theories for calculation the probability \eqref{EQ1} have been developed, most of them are only concentrated on the constant boundary, and even in this case, few of these methods got an brief explicit analytic formulas. For studying the details of these results, we refer the reader to \cite{jamison1970reciprocal},\cite{shepp1971first},\cite{shepp1971first},\cite{shepp1976first} and \cite{orsingher1989maximum} and the references therein. Analytic formulas of probability \eqref{EQ1} associated with linear and piecewise linear boundary function $f$ was developed by Ein-Gal and Bar-David, Abrahams, Bischoff (see \cite{10.2307/2959476}, \cite{abrahams1984ramp} and \cite{Bischoff2016} respectively). All those formulas  are based on the Markov-like property (or reciprocal property) of $S$ stated by Jamison in \cite{jamison1970reciprocal}.
\newline
In this paper, we use a new method to derive an explicit formula for the non-crossing probabilities of Slepian  process $S$ with piecewise linear boundary $f$. This formula is used to develop a general formula for general continuous boundaries by approximating the boundaries with piecewise linear function. The formula we develop for piecewise linear boundary is equivalent to the results in \cite{Bischoff2016}, but has a simper form and is tractable for applications.

The paper is organized as follows: In Section \ref{Preliminaries}, we give some preliminaries which we will use to derive calculation the non-crossing probabilities of Slepian  process, In Section \ref{Somere}, we use the results in Section \ref{Preliminaries} to derive explicit formulas for the non-crossing probabilities of Slepian, several cases are involved. Finally, proofs of lemmas are given in Section \ref{Proofs}.

\section{ Preliminaries} \label{Preliminaries}
Throughout this paper, we consider the one-dimensional Slepian process defined as
the increment of a Brownian motion process, namely
\BQN\label{EQ2.1}
S(t)=B(t+1)-B(t),\quad t\in [0,1],
\EQN
where $B(t)$ is a standard Brownian motion. It is easy to verify that $S(t),t\in [0,1]$ is a centered stationary Gaussian process with covariance function
\BQNY
R_{S}(t_1,t_2)=\E({S(t_1)S(t_2)})=1-\abs{t_1-t_2},\quad t_1,t_2\in [0,1].
\EQNY
For any general boundary function $f$, the direct calculation of the non-crossing probability is always impossible. One feasible method is to derive some tractable formula for $p_f$ with $f$ a a piecewise linear function. Then
the general boundary case can be determined by approximating  $f$ by piecewise linear functions, see for instance
\cite{wang1997boundary,1137.60023, 1079.62047, 1103.60040}.

However, unlike the Brownian motion, calculating the non-crossing probability for Slepian processes is more complicated because Brownian motion is Markovian, but the Slepian process is not (recall that the only non-trivial Gaussian Markov process is the Ornstein-Uhlenbeck process, see e.g. \cite{doob1942brownian}). Nonetheless, in this contribution we first derive a  direct relationship between Brownian motion and Slepian processes which is useful for calculation of $p_f$.

\BL\label{lemma2.1}
If $S=\{S(t),\;t\in [0,1]\}$ is a Slepian process, then for any $t\in [0,1]$ the conditional density of $S(t)$ given $S(0)$ is given by
\BQN\label{EQ2.2}
\varphi(S(t)=y\mid S(0)=x)=\frac{1}{\sqrt{2\pi t(2-t)}}\exp \Bigr\{-\frac{1}{2t(2-t)}(y+(t-1)x)^2\Bigl\}, \quad x,y\in \R.
\EQN
\EL
The proof of this lemma is presented in the Section \ref{Proofs}. From Lemma \ref{lemma2.1} the  conditional Slepian process obtained by conditioning on $\{ S(0)=x\}$ is also a Gaussian process. Concretely, if we denote the process $Y=\Bigr\{Y_t=(S(t)\mid S(0)=x),\;t\in [0,1]\Bigl\}$, then we have the following corollary:

\BK\label{korr2.2}
$Y=\Bigr\{Y_t,\;t\in [0,1]\Bigl\}$ is a Gaussian process with mean $\E(Y_t)=(1-t)x$ and covariance function
\BQNY
R_{Y}(t_1,t_2)=\E({Y_{t_1}Y_{t_2}})=\min\{t_1(2-t_2),t_2(2-t_1)\}.
\EQNY
\EK
The proof of the above claim can be directly shown by calculating the covariance function of $Y$.   The covariance function of $Y$ has a minimum form which is the case in Brownian motion, hence we can construct a new unconditional process from Brownian motion which has the same covariance function with process $Y$.
\BL\label{lemma2.3}
Suppose $B=\Bigr\{B_t,\;t\in [0,1]\Bigl\}$ is a standard Brownian motion, let $Z=\Bigr\{Z_t=(2-t)B(\frac{t}{2-t})+(1-t)x,\;t\in [0,1]\Bigl\}$. Then $Z$ is a Gaussian process with mean $\E(Z_t)=(1-t)x$ and covariance function
\BQNY
R_{Z}(t_1,t_2)=\E({Z_{t_1}Z_{t_2}})=\min\{t_1(2-t_2),t_2(2-t_1)\}.
\EQNY
\EL
We present the proof of this lemma in Section \ref{Proofs}. Combining Corollary \ref{korr2.2} and Lemma \ref{lemma2.3}, we obtain that $Y$ is equivalent in distribution with $Z$ i.e. $Y\equivdis Z$.
\BRM
The processes $Z$ is just a space-time transformation of standard Brownian, it is thus more flexible due to the rich available results in Brownian motion.
\ERM
\section{ Boundary Non-Crossing Probability} \label{Somere}
We have obtained that under condition $S(0)=x$, the conditional distribution of Slepian process $S(t)$ is equivalent in distribution with a space-time transformation of Brownian motion. In this section, we use this fact to develop the boundary non-crossing probabilities of Slepian process. Four cases will be involved depend on the boundary function $f$, namely, constant boundary, linear boundary, piecewise linear boundary and general boundary. Though the piecewise linear boundary include the former two cases, we consider them separately for the reason that the non-crossing probabilities in the former two cases have a formula which can be represented by the density and cumulative function of standard normal distribution, and these formulas are more simper compare with the corresponding results in the literature.
\subsection{Constant boundary} \label{constant}
The constant boundary is the simplest case, we first cite a well known result obtained by Shepp in \cite{shepp1971first} (see also in \cite{cressie1980asymptotic}). Suppose $a\in\R$, then
\BQNY
\pk{S(t)\leq a,\;\text{for all}\;t\in [0,1]}&=&\int_{D'}\det\left (\begin{array}{cc}
                                                               \phi(a-y_1) & \phi(a-y_2) \\
                                                               \phi(a) & \phi(y_1-y_2+a)
                                                             \end{array}\right )dy_1dy_2\\
&=&\int_{D'}\phi(a-y_1)\phi(y_1-y_2+a)-\phi(a)\phi(a-y_2)dy_1dy_2,
\EQNY
where $\phi$ is the density function of a standard normal distribution $N(0,1)$ and $D'=\{0<y_1<y_2\}$. This formula is suitable for numerical calculation. However, the integral over the  region $D'$ is not explicite;  our nexst result gives a tractable simpler formula:
\BT\label{T3.1}
For $a \in \R$, the non-crossing probabilities of \eqref{EQ1} is given by
\BQN\label{eq3.0}
\pk{S(t) \leq a,\;\text{for all}\;t\in [0,1]}=\Phi^2(a)-a\phi(a)\Phi(a)-\phi^2(a),
\EQN
where $\phi$ is the density function of standard normal distribution, i.e. $\phi(x)=\frac{1}{\sqrt{2\pi}}e^{-\frac{x^2}{2}}$, and $\Phi(x)=\int_{-\infty}^{x}\phi(s)ds$ is the cumulative distribution function of standard normal distribution.
\ET
The proof of this theorem is displayed in Section \ref{Proofs}.
\BRM\label{remark1}
In Theorem \ref{T3.1}, if $a=0$, then the probabilities that $S(t)\leq 0, t\in [0,1]$ is
\BQN
\pk{\sup_{0\leq t\leq 1}S(t) \leq 0}=\Phi^2(0)-\phi^2(0)=\frac{1}{4}-\frac{1}{2\pi}.
\EQN
Furthermore, if we define the stoping time $\tau=\inf\{t\geq 0,S(t)=0\}$, then
\BQN
\pk{\tau\leq 1}=1-2(\frac{1}{4}-\frac{1}{2\pi})=\frac{1}{\pi}-\frac{1}{2}.
\EQN
\ERM
\BRM
For $a\geq 0$, from Theorem \ref{T3.1} we have the inequality that
\BQNY
\pk{\sup_{0\leq t\leq 1}\abs{S(t)} \leq a}<\pk{\abs{\sup_{0\leq t\leq 1}S(t)} \leq a}=2\Phi(a)-a\phi(a)-1
\EQNY
\ERM
\subsection{Linear boundary} \label{linear}
Various results for calculating the non-crossing probability with linear boundary are only about the Brownian motion, among them one famous is Bachelier-Levy formula (see e.g. \cite{deelstra1994remarks}) which is given by the following,
for $a>0$
\BQN\label{EQ4.1}
\pk{B(t)\leq a+bt,\;\text{for all}\;t\in [0,T]}=\Phi(b\sqrt{T}+\frac{a}{\sqrt{T}})-e^{-2ab}\Phi(b\sqrt{T}-\frac{a}{\sqrt{T}})
\EQN
and the above probability is 0 when $a\leq 0$. Since the proof of Equation \eqref{EQ4.1} used the Markovian property of Brownian motion, thus can not be applied to Slepian process. Fortunately, we have the following similar result.
\BT\label{T3.4}
For any $a\in \R,\;b\neq0$, the non-crossing probabilities of \eqref{EQ1} is given by
\BQN\label{EQ3.4}
\pk{S(t) \leq a+bt,\;\text{for all}\;t\in [0,1]}=\Phi(a+b)\Phi(a)-\frac{1}{b}\phi(a)\Phi(a+b)+\frac{\sqrt{2\pi}}{b}\phi(a)\phi(b)\Phi(a)e^{-ab},
\EQN
where $\phi$ is the density function of standard normal distribution, i.e. $\phi(x)=\frac{1}{\sqrt{2\pi}}e^{-\frac{x^2}{2}}$, and $\Phi(x)=\int_{-\infty}^{x}\phi(s)ds$ is the cumulative distribution function of standard normal distribution.
\ET
The proof of this theorem is given in section \ref{Proofs}.
\BRM
(i) Particularly in Theorem \ref{T3.4}, if $a=0,b\neq 0$, then the probabilities that $S(t)\leq bt, t\in [0,1]$ is
\BQNY
\pk{\sup_{0\leq t\leq 1}S(t) \leq bt}=(\frac{1}{2}-\frac{1}{b})\Phi(b)+\frac{\sqrt{2\pi}}{2b}\phi(b).
\EQNY
(ii) The result in Theorem \ref{T3.4} in the limiting case as $b\to 0$ coincide with the result in constant case.
\ERM
\subsection{Piecewise linear boundary} \label{piecewise}
We now generalize the formula \eqref{EQ3.4} to a piecewise linear boundary, in this subsection, the boundary function $f$ is piecewise, to be specific, let $0=t_0<t_1<\ldots<t_{n-1}<t_n=1$ be a partition of $[0,1]$, and $f(t)$ be a continuous, and linear function on each of interval $[t_{i-1},t_i],\;i=1,2,\ldots,n$. The method we use based on an classical result proved by Wang and Potzelberger, when calculating non-crossing probability of standard Brownian motion, which we cite without proof as following
\BL\label{Wang}
\cite{wang1997boundary}The non-crossing probability of standard Brownian motion with boundary $f(t)$ is given by
\BQN\label{EQ3.1}
&&\pk{B_u \leq c(u),\;\text{for all}\;u\in [0,1]}\\
\nonumber &=&\int_{-\infty}^{c_1}\int_{-\infty}^{c_2}\ldots\int_{-\infty}^{c_n}
(2\pi)^{-\frac{n}{2}}\Pi_{i=1}^{n}\frac{1}{\sqrt{t_i-t_{i-1}}}\exp\Bigr\{-\frac{(u_i-u_{i-1})^2}{2(t_i-t_{i-1})}\Bigl\}\\
\nonumber&\times&\Pi_{i=1}^{n}(1-\exp\Bigr\{-\frac{2(u_{i-1}-c_{i-1})(u_i-c_{i})}{t_i-t_{i-1}}\Bigl\})du_ndu_{n-1}\ldots du_1,
\EQN
where $c_i:=c(t_i),i=0,1,2,\ldots,n$ for short.
\EL
Noting that under the space-time transformation $g(t)=(2-t)f(\frac{t}{2-t})+(1-t)x$, the function $g$ is also a piecewise linear function on $[0,1]$, hence the non-crossing probability of Slepian process can transform into the non-crossing probability of Brownian motion, in fact, we have
\BT\label{T3.5}
Let $0=t_0<t_1<\ldots<t_{n-1}<t_n=1$ and $l(t)$ be continuous function and linear on each of interval $[t_{i-1},t_i],\;i=1,2,\ldots,n$.
Suppose the parameters of line segment on interval $[t_{i-1},t_i]$ are $a_i,b_i$. Then the non-crossing probability of Slepian process for piecewise linear boundary of $l(t)$ is given by
\BQNY
&&\pk{S(t) \leq l(t),\;\text{for all}\;t\in [0,1]}\\
&=&\int_{-\infty}^{a_1}\int_{-\infty}^{h_1(x)}\int_{-\infty}^{h_2(x)}\ldots\int_{-\infty}^{h_n(x)}
(2\pi)^{-\frac{n+1}{2}}\exp\Bigr\{-\frac{1}{2}x_{0}^2\Bigl\}\Pi_{i=1}^{n}\frac{1}{\sqrt{t_i-t_{i-1}}}\exp\Bigr\{-\frac{(u_i-u_{i-1})^2}{2(t_i-t_{i-1})}\Bigl\}\\
&\times&\Pi_{i=1}^{n}(1-\exp\Bigr\{-\frac{2(u_{i-1}-h_{i-1}(x))(u_i-h_{i}(x))}{t_i-t_{i-1}}\Bigl\})du_n du_{n-1}\ldots du_1dx,
\EQNY
where $h_{i}(x)=(\frac{a_i}{2}+b_i+\frac{x}{2})t_i+\frac{a_i}{2}-\frac{x}{2},\; i=1,2,\ldots n$, and $h_{0}(x)=\frac{a_1-x}{2},\;u_0=0$.
\ET
The proof of this theorem is displayed in section \ref{Proofs}.
\subsection{General boundary}\label{general}
Observing that any continuous functions on $[0,1]$ can be approximated by continuous piecewise linear function uniformly, then we can develop the formula for the non-boundary probability of Slepian process with general boundary in a limiting form. In fact, we have the following
\BT\label{T3.6}
Suppose $f(u)$ be a continuous functions on $[0,1]$, and $l_{n}(u)$ defined on $[0,1]$ are continuous piecewise linear function which admire $l_{n}(u)\to f_{u}$ uniformly as $n\to \infty$. Then
\BQNY
&&\pk{S(t) \leq f(t),\;\text{for all}\;t\in [0,1]}\\
&=&\lim_{n\to\infty}\int_{-\infty}^{a_0}\int_{-\infty}^{h_1(x)}\int_{-\infty}^{h_2(x)}\ldots\int_{-\infty}^{h_n(x)}
(2\pi)^{-\frac{n+1}{2}}\exp\Bigr\{-\frac{1}{2}x_{0}^2\Bigl\}\Pi_{i=1}^{n}\frac{1}{\sqrt{t_i-t_{i-1}}}\exp\Bigr\{-\frac{(u_i-u_{i-1})^2}{2(t_i-t_{i-1})}\Bigl\}\\
&\times&\Pi_{i=1}^{n}(1-\exp\Bigr\{-\frac{2(u_{i-1}-h_{i-1}(x))(u_i-h_{i}(x))}{t_i-t_{i-1}}\Bigl\})du_ndu_{n-1}\ldots du_1dx,
\EQNY
where $h_{i}(x)=(\frac{a_i}{2}+b_i+\frac{x}{2})t_i+\frac{a_i}{2}-\frac{x}{2},\; i=1,2,\ldots n$, and $h_{0}(x)=\frac{a_1-x}{2},\;u_0=0$.
\ET
The proof of this theorem is given in Section \ref{Proofs}.
\subsection{An equivalent result}\label{compare}
It is worth mentioning that Theorem \ref{T3.6} is a result for general boundary function, the spacial case when $g(u)$ is piecewise function stated in Theorem \ref{T3.5} which is equivalent to Theorem 2.1 proved by Bischoff in \cite{Bischoff2016}. In this subsection we will talk about this equivalence and state that our formulas is convenience in calculation, precisely, for any $h\leq 1$ and continuous boundary function $g$, applying Fubini's Theorem we have
\BQN\label{equi}
&&\int_{-\infty}^{g(0)}\int_{-\infty}^{g(h)}\pk{S(t)\leq g(t),\;\text{for all}\;t\in [0,h]\mid S(0)=x,S(h)=x_h}\varphi(S(0)=x,S(h)=x_h)dx_h dx\\
\nonumber &=&\pk{S(t)\leq g(t),\;\text{for all}\;t\in [0,h]}\\
\nonumber &=&\int_{-\infty}^{g(0)}\pk{S(t)\leq g(t),\;\text{for all}\;t\in [0,h]\mid S(0)=x}
\varphi(S(0)=x)dx.
\EQN
When $g(t)$ is a piecewise linear function, the first equality above derives Theorem 2.1 in \cite{Bischoff2016} while the second equality derives Theorem \ref{T3.5} in our text. Furthermore, our result in Theorem \ref{T3.5} is more convenient than Theorem 2.1 and Theorem 2.2 in \cite{Bischoff2016} in calculating the boundary non-crossing probabilities. To see this, we compute the probability that $S(t)$ not crossing the boundary zero for all $u\in [0,1]$ with Bischoff's method. In fact, we can rewrite the probability that
\BQNY
&&\pk{S(t)\leq 0,\;\text{for all}\;t\in [0,1]} \\
&=&\int_{-\infty}^{0}\int_{-\infty}^{0}\pk{S(t)\leq 0,\;\text{for all}\;t\in (0,1)\mid S(0)=x,S(1)=x_1}\varphi(S(0)=x,S(1)=x_1)dx_1dx,
\EQNY
where $\varphi(S(0)=x,S(1)=x_1)=\frac{1}{2\pi}\exp\Bigr\{-\frac{1}{2}(x_{1}^{2}+x_{0}^{2})\Bigl\}$ is the density of the joint distribution of $(S(0),S(1))$.
\newline
From the Theorem 2.1 and Theorem 2.2 in \cite{Bischoff2016} (also see in \cite{10.2307/2959476}), we obtain
\BQNY
&&\pk{S(t)\leq 0,\;\text{for all}\;t\in (0,1)\mid S(0)=x,S(1)=x_1}\\
&=&1+\frac{x}{2\sqrt{\pi}}\exp\Bigr\{\frac{(x_1-x)^2}{4}\Bigl\}\int_{0}^{1}\frac{1}{\sqrt{(1-s)s^3}}\exp\Bigr\{-\frac{1}{4}(\frac{x_{0}^{2}}{s}+\frac{x_{1}^{2}}{1-s})\Bigl\}ds.
\EQNY
Letting $\sigma^2=\frac{2(1-s)}{2-s}$, therefore,
\BQNY
&&\pk{S(t)\leq 0,\;\text{for all}\;t\in [0,1]}\\
&=&\frac{1}{4}+\frac{1}{2\pi}\int_{-\infty}^{0}\int_{-\infty}^{0}\frac{1}{2\pi}\exp\Bigr\{-\frac{1}{2}(x_{1}^{2}+x_{0}^{2})\Bigl\}\frac{x}{2\sqrt{\pi}}\exp\Bigr\{\frac{(x_1-x)^2}{4}\Bigl\}
\int_{0}^{1}\frac{1}{\sqrt{(1-s)s^3}}\exp\Bigr\{-\frac{1}{4}(\frac{x_{0}^{2}}{s}+\frac{x_{1}^{2}}{1-s})\Bigl\}dsdx_1dx\\
&=&\frac{1}{4}+\frac{1}{2\pi}\int_{-\infty}^{0}\int_{-\infty}^{0}\int_{0}^{1}\frac{x}{2\sqrt{\pi}}\exp\Bigr\{-\frac{x_{0}^{2}}{2s(2-s)}\Bigl\}\frac{1}{\sqrt{(1-s)s^3}}
\exp\Bigr\{-\frac{2-s}{4(1-s)}(x_1+\frac{1-s}{2-s}x)^2\Bigl\}dsdx_1dx\\
&=&\frac{1}{4}+\frac{1}{2\pi}\int_{0}^{1}\int_{-\infty}^{0}x\frac{1}{\sqrt{(1-s)s^3}}\Phi(\frac{1-s}{2-s}\frac{x}{\sigma})\exp\Bigr\{-\frac{x_{0}^{2}}{2s(2-s)}\Bigl\}dxds\\
&=&\frac{1}{4}-\int_{0}^{1}\int_{-\infty}^{0}\frac{1}{2\pi}\sqrt{\frac{2-s}{s}}\Phi(\frac{1-s}{2-s}\frac{x}{\sigma})de^{-\frac{x_{0}^{2}}{2s(2-s)}}ds\\
&=&\frac{1}{4}-\frac{1}{4\pi}\int_{0}^{1}\sqrt{\frac{2-s}{s}}ds+\frac{1}{4\pi}\int_{0}^{1}\sqrt{\frac{1-s}{1+s}}ds\\
&=&\frac{1}{4}-\frac{1}{4\pi}(1+\frac{\pi}{2})+\frac{1}{4\pi}(\frac{\pi}{2}+1)\\
&=&\frac{1}{4}-\frac{1}{2\pi}.
\EQNY
However, from Remark \ref{remark1}, this is obvious in our text.
\newline
In the following, we let $g(t)=a+bt,t\in [0,h]$ be a linear boundary function. We let the first hitting time $\tau_g=\inf\{t\geq0:S(t)>g(t)\}$ and denote $\pi_g(\cdot\mid S(0)=x)$ the Lebesgue-density of the first hitting time under the condition that the process $S$ starts in $x$ at $t=0$. For any $0\leq s_1<s_2<\ldots<s_m\leq 1, 0\leq t_1<t_2<\ldots<t_n\leq 1, s_i\neq t_j$, denote $\varphi(S(s_1),s(s_2),\cdots,s(s_m))$ the density of the finite dimensional distribution of $(S(s_1),s(s_2),\cdots,s(s_m))$ and $\varphi(S(s_1),s(s_2),\cdots,s(s_m)\mid S(t_1),s(t_2),\cdots,s(t_n))$ the conditional finite dimensional distribution of $(S(s_1),s(s_2),\cdots,s(s_m))$ given $S(t_1),s(t_2),\cdots,s(t_n)$.
Then we have the following corollary
\BK\label{corollary2}
With all notations above, we have the relationship that
\BQNY
&&\int_{-\infty}^{a}\int_{-\infty}^{a+bh}\varphi(S(0)=x,S(h)=x_h)dx_hdx-\int_{-\infty}^{a}\int_{-\infty}^{a+bh}\int_{0}^{h}\frac{a-x}{u}\varphi(S(u)=a+bu, S(0)=x,S(h)=x_h)dudx_hdx\\
\nonumber&=&\int_{-\infty}^{a}\varphi(S(0)=x)dx-\int_{-\infty}^{a}\int_{0}^{h}\frac{a-x}{u}\varphi(S(u)=a+bu, X=x)dudx.
\EQNY
\EK
\section{Proofs}\label{Proofs}
\prooflem{lemma2.1}
Since $S(t)$ is a  stationary Gaussian process, then for any $t\in [0,1]$, from equation (3) in \cite{slepian1961}, the joint density of $(S(t),S(0))$ was given by
\BQNY
\varphi(S(t)=y, S(0)=x)
=\frac{1}{2\pi\sqrt{t(2-t)}}\exp\Bigr\{-\frac{1}{4}\Bigr(\frac{(y+x)^2}{2-t}+\frac{(y-x)^2}{t}\Bigl)\Bigl\}.
\EQNY
Hence the conditional density of $S(t)$ given $S(0)$ is
\BQNY
\varphi(S(t)=y\mid S(0)=x)&=&\frac{\varphi(S(t)=y,S(0)=x)}{\varphi(S(0)=x)}\\
&=&\frac{1}{\sqrt{2\pi t(2-t)}}\exp \Bigr\{-\frac{1}{2t(2-t)}(y+(t-1)x)^2\Bigl\}.
\EQNY
\QED

\prooflem{lemma2.3}
The proof of $Z_t,\;t\in [0,1]$ is a Gaussian process is obvious, and followed from the property of Standard Brownian motion $B$. Without loss of generosity, suppose $0\leq t_1\leq t_2\leq 1$, then
\BQNY
R_{Z}(t_1,t_2)&:=&\E[Z_{t_1}Z_{t_2}]-\E{Z_{t_1}}\E{Z_{t_2}}\\
&=&\E\Bigl\{[(2-t_1)B(\frac{t_1}{2-t_1})+(1-t_1)x][(2-t_2)B(\frac{t_2}{2-t_2})+(1-t_2)x]\Bigr\}-(1-t_1)(1-t_2)x_{0}^{2}\\
&=&(2-t_1)(2-t_2)\cdot\min(\frac{t_1}{2-t_1},\frac{t_2}{2-t_2})\\
&=&t_1(2-t_2).
\EQNY
Completing the proof.
\QED

\prooftheo{T3.1}
Observing that the probability \eqref{EQ1} can be rewritten as
\BQNY
\pk{S(t)\leq a,\;\text{for all}\;t\in [0,1]}&=&\int_{-\infty}^{a}\pk{S(t)\leq a,\;\text{for all}\;t\in [0,1]\mid S(0)=x}
\varphi(S(0)=x)dx
\EQNY
and the conditional distribution of $S(t)$ given $S(0)$ is equivalent to the distribution of process $Z(t)$, then
\BQNY
\pk{S(t)\leq a,\;\text{for all}\;t\in [0,1]}
&=&\int_{-\infty}^{a}\pk{Z_u\leq a,\;\text{for all}\;u\in [0,1]}\phi(x)dx,
\EQNY
where $\phi$ is the density function of standard normal distribution, i.e. $\phi(x)=\frac{1}{\sqrt{2\pi}}e^{-\frac{x^2}{2}}$, and denote $\Phi(x)=\int_{-\infty}^{x}\phi(s)ds$ the cumulative distribution function of standard normal distribution. Then
\BQNY
\pk{S(t)\leq a,\;\text{for all}\;t\in [0,1]}
&=&\int_{-\infty}^{a}\P[(2-t)B(\frac{t}{2-t})+(1-t)x\leq a,\;\text{for all}\;t\in [0,1]]d\Phi(x)\\
&=&\int_{-\infty}^{a}\P[B(u)\leq(\frac{a+x}{2})u+\frac{a-x}{2} ,\;\text{for all}\;u\in [0,1]]d\Phi(x).
\EQNY
Since, the probability that $B(t)\leq a+bt$ for all $t\in [0,T]$ is well known, and can be obtained from the famous Bachelier-Levy formula (see Equation \eqref{EQ4.1}). Hence, the probability that $S(t)\leq a$ for all $t\in [0,1]$ can be written as following
\BQNY
\pk{S(t)\leq a,\;\text{for all}\;t\in [0,1]}
&=&\int_{-\infty}^{a}\Bigl\{\Phi(a)-e^{-\frac{a^2-x_{0}^{2}}{2}}\Phi(x)\Bigr\}d\Phi(x)\\
&=&\Phi^2(a)-e^{-\frac{a^2}{2}}\int_{-\infty}^{a}e^{-\frac{x_{0}^{2}}{2}}\Phi(x)d\Phi(x)\\
&=&\Phi^2(a)-\phi(a)\int_{-\infty}^{a}\Phi(x)dx\\
&=&\Phi^2(a)-a\phi(a)\Phi(a)-\phi^2(a),
\EQNY
which completes the proof.
\QED

\prooftheo{T3.4}
The proof of this theorem is similar to Theorem \ref{T3.1}. Using the same definitions of $\phi$ and $\Phi$ as before, we compute directly the probability that $S(t)\leq a$ for all $t\in [0,1]$.
\BQNY
\pk{S(t)\leq a+bt,\;\text{for all}\;t\in [0,1]}&=&\int_{-\infty}^{a}\pk{S(t)\leq a+bt,\;\text{for all}\;t\in [0,1]\mid S(0)=x}
\varphi(S(0)=x)dx\\
&=&\int_{-\infty}^{a}\P[Z_u\leq a+bu,\;\text{for all}\;u\in [0,1]]\phi(x)dx.
\EQNY
Substituting $Z_u=(2-u)B(\frac{u}{2-u})+(1-u)x$ into the above and after some linear transformation we have
\BQNY
\pk{S(t)\leq a+bt,\;\text{for all}\;t\in [0,1]}&=&\int_{-\infty}^{a}\pk{B(t)\leq \frac{a-x}{2}+(\frac{a+x}{2}+b)t,\;\text{for all}\;t\in [0,1]}d\Phi(x).
\EQNY
Using Bachelier-Levy formula again, we therefore obtain
\BQNY
\pk{S(t)\leq a+bt,\;\text{for all}\;t\in [0,1]}&=&\int_{-\infty}^{a}\Phi(a+b)-\exp\{-\frac{a^2-x_{0}^{2}}{2}-ab+bx\}\Phi(b+x)d\Phi(x)\\
&=&\Phi(a+b)\Phi(a)-\frac{1}{\sqrt{2\pi}}\exp\{-\frac{a^2}{2}-ab\}\int_{-\infty}^{a}\Phi(b+x)\exp\{bx\}dx.
\EQNY
Let $I(a,b)=\int_{-\infty}^{a}\Phi(b+x)\exp\{bx\}dx$, observing that $b\neq 0$, then
\BQNY
I(a,b)&=&\frac{1}{b}\int_{-\infty}^{a}\Phi(b+x)d\exp\{bx\}\\
&=&\frac{1}{b}\Bigl\{\Phi(a+b)e^{ab}-\frac{1}{\sqrt{2\pi}}\int_{-\infty}^{a}\exp\{-\frac{(b+x)^2}{2}\}\exp\{bx\}\Bigr\}dx\\
&=&\frac{1}{b}\Bigl\{\Phi(a+b)e^{ab}-\exp\{-\frac{b^2}{2}\}\Phi(a)\Bigr\}\\
&=&\frac{1}{b}\Bigl\{\Phi(a+b)e^{ab}-\sqrt{2\pi}\phi(b)\Phi(a)\Bigr\}.
\EQNY
Thus,
\BQNY
\pk{S(t)\leq a+bt,\;\text{for all}\;t\in [0,1]}=\Phi(a+b)\Phi(a)-\frac{1}{b}\phi(a)\Phi(a+b)+\frac{\sqrt{2\pi}}{b}\phi(a)\phi(b)\Phi(a)e^{-ab},
\EQNY
hence the proof is complete. \QED

\prooftheo{T3.5}
With the result of linear boundary case in mind, (see Theorem \ref{T3.4}), the proof for this theorem is easy. In fact, with the same method above, we can rewrite the probability that $S(t)\leq l(t)$ for all $t\in [0,1]$ as the integral of  some Brownian non-crossing probability. Concretely, we have
\BQNY
\pk{S(t)\leq l(t),\;\text{for all}\;t\in [0,1]}&=&\int_{-\infty}^{l(0)}\P[S(t)\leq l(t),\;\text{for all}\;t\in [0,1]\mid S(0)=x]
\varphi(S(0)=x)dx\\
&=&\int_{-\infty}^{l(0)}\P[B(t)\leq \frac{t+1}{2}l(\frac{2t}{t+1})-\frac{1-t}{2}x,\;\text{for all}\;t\in [0,1]]d\Phi(x).
\EQNY
Since $l(t)$ be continuous function and piecewise linear with parameters $a_i,b_i$ on each of interval $[t_{i-1},t_i],\;i=1,2,\ldots,n$, i.e. $l(t)=a_i+b_it,\;t\in [t_{i-1},t_i],i=1,2,\ldots, n.$ where $0=t_0<t_1<\ldots<t_{n-1}<t_n=1$, We have $l(0)=a_1$.\\
Next, given $x$, we define $h(x,t)=\frac{t+1}{2}l(\frac{2t}{t+1})-\frac{1-t}{2}x$, then $h(x,t)$ is also a piecewise linear function about $t$, and have the following form
\BQNY
h(x,t)=(\frac{a_i+x}{2}+b_i)t+\frac{a_i-x}{2},\quad t\in [t_{i-1},t_i],i=1,2,\ldots,n.
\EQNY
Denoting $h_{i}(x):=h(x,t_i)=(\frac{a_i}{2}+b_i+\frac{x}{2})t_i+\frac{a_i}{2}-\frac{x}{2},\; i=1,2,\ldots n$, and $h_{0}(x)=\frac{a_1-x}{2},u_0=0$ for short, By Lemma \ref{Wang}, we have
\BQNY
&&\pk{S(t) \leq l(t),\;\text{for all}\;t\in [0,1]}\\
&=&\int_{-\infty}^{a_1}\int_{-\infty}^{h_1(x)}\int_{-\infty}^{h_2(x)}\ldots\int_{-\infty}^{h_n(x)}
(2\pi)^{-\frac{n+1}{2}}\exp\Bigr\{-\frac{1}{2}x_{0}^2\Bigl\}\Pi_{i=1}^{n}\frac{1}{\sqrt{t_i-t_{i-1}}}\exp\Bigr\{-\frac{(u_i-u_{i-1})^2}{2(t_i-t_{i-1})}\Bigl\}\\
&\times&\Pi_{i=1}^{n}(1-\exp\Bigr\{-\frac{2(u_{i-1}-h_{i-1}(x))(u_i-h_{i}(x))}{t_i-t_{i-1}}\Bigl\})du_ndu_{n-1}\ldots du_1dx,
\EQNY
establishing the proof.
\QED

\prooftheo{T3.6}
Combining Theorem \ref{T3.5}, this is an obvious consequence from the continuity of probability measure $\P$ and the uniform convergence of $l_{n}(u)$ to $f(u)$.
\QED

\proofkorr{corollary2}
From (A.4) in \cite{Bischoff2016}, we have
\BQNY
\pi_g(t\mid S(0)=x)&=&\frac{a-x}{t}\varphi(S(t)=a+bt\mid S(0)=x),\quad t\in (0,h],\\
\pi_g(t\mid S(0)=x,S(h)=x_h)&=&\frac{\varphi(S(h)=x_h\mid S(0)=x,S(t)=a+bt)\cdot\pi_g(t\mid S(0)=x)}{\varphi(S(h)=x_h\mid S(0)=x)}.
\EQNY
We then can represent the probability that $S(t)\leq a+bt$ for all $u\in [0,h]$ given $S(0)=x$ and $S(h)=x_h$ as
\BQNY
\pk{S(t)\leq a+bt,\;\text{for all}\;t\in [0,h]\mid S(0)=x,S(h)=x_h}=1-\int_{0}^{h}\frac{a-x}{u}\varphi(S(u)=a+bu\mid S(0)=x,S(h)=x_h)du.
\EQNY
Similarly, the probability that $S(t)\leq a+bt$ for all $t\in [0,h]$ given $S(0)=x$ can be written as
\BQNY
\pk{S(t)\leq a+bt,\;\text{for all}\;t\in [0,h]\mid S(0)=x}&=&1-\int_{0}^{h}\frac{a-x}{u}\varphi(S(u)=a+bu\mid S(0)=x)du,
\EQNY
the conclusion then follows from Equation \eqref{equi}.
\QED
\section{Acknowledgement}
This work was partly financed by the project NSFC Grant NO.71573143 and SNSF Grant 200021-166274.
\bibliographystyle{ieeetr}
\bibliography{Slep4}
\end{document}